\documentclass{amsart}

\usepackage{amsmath, amssymb, amsthm, mathtools}

\newtheorem{theorem}{Theorem}[section]
\newtheorem{proposition}[theorem]{Proposition}
\newtheorem{lemma}[theorem]{Lemma}

\theoremstyle{definition}
\newtheorem{example}[theorem]{Example}
\newtheorem{definition}[theorem]{Definition}

\theoremstyle{remark}
\newtheorem{remark}[theorem]{Remark}

\title[Ad-conditions for higher-order recurrences]{Ad-conditions for higher-order recurrences: Dunkl and matrix-valued settings}
\author{I. Bono Parisi}
\address{CIEM-FaMAF, Universidad Nacional de C\'ordoba, CP 5000, C\'ordoba, Argentina}
\email{ignacio.bono@unc.edu.ar}
\author{F. A. Gr\"unbaum}
\address{Department of Mathematics, University of California, Berkeley, CA 94720, USA}
\email{grunbaum@math.berkeley.edu}
\author{I. Zurri\'an}
\address{Departamento de Matem\'atica Aplicada II and IMUS, Universidad de Sevilla, Sevilla, Spain}
\email{ignacio.zurrian@fulbrightmail.org}

\date{}

\begin{document}

\begin{abstract}
The ad-conditions of Duistermaat and Gr\"unbaum characterize bispectrality for a pair of differential operators. M.~Reach observed that a sequence of functions which is simultaneously an eigenfunction of a second-order differential operator and satisfies a recurrence relation of arbitrary order must satisfy a different kind of ad-conditions. We prove that the Reach ad-conditions are also sufficient, by an argument which is purely algebraic. This motivated us to explore ad-conditions in two new directions: to Dunkl differential-difference operators, where the coefficients of the ad-condition become elements of the Dunkl algebra and the reflection enters as a coefficient, and to matrix-valued orthogonal polynomials, where we exhibit ad-conditions for several families of exceptional matrix-valued polynomials satisfying five- and seven-term recurrence relations.
\end{abstract}

\maketitle

\section{Introduction}

We say that we have a bispectral situation when a family of functions
$\psi(x,k)$ consists of eigenfunctions of two different operators, one
acting in $x$ and the other one in $k$; in each case the operator depends
on only one of the variables and the eigenvalue on the other one.

The bispectral problem for functions $\psi(x,k)$ when both variables are
continuous was considered in \cite{DG}. This was the first appearance of
the so-called \emph{ad-conditions} as a necessary and sufficient condition
for bispectrality, involving a pair of differential operators, one of order
two and the other one of arbitrary order. In \cite{GH96,GH} one sees how the use
of these conditions leads to the ``classical orthogonal polynomials'',
whose study goes back at least to S.~Bochner \cite{Boch}; here a
differential and a difference operator, both of order two, are involved.
In \cite{R1,R2}, based on his Berkeley thesis, M.~Reach deals with general bispectral sequences $p_{n}(x)$,
with a differential operator of order two and a difference operator of
arbitrary order. He observes that they necessarily satisfy a different kind
of ad-conditions, which we call \emph{the Reach ad-conditions}.

In this paper we do several things. In Section~\ref{sec:scalar} we prove
that, in the scalar case, the Reach ad-conditions are also sufficient for
bispectrality (Theorems~\ref{thm:three-term} and~\ref{thm:BGZ-general}).
The proof is purely algebraic, and this is what makes the rest of the
paper possible. In Section~\ref{sec:dunkl} we show that the conditions
survive when the differential operator is replaced by a Dunkl
differential-difference operator \cite{Dun89}; the new phenomenon is that the
coefficients of the ad-condition are no longer scalars but elements of the
Dunkl algebra, with the reflection appearing as a coefficient. In
Section~\ref{sec:matrix} we consider matrix-valued orthogonal polynomials,
where both the differential and the difference operators have matrix-valued
coefficients, and we exhibit ad-conditions for several families of
exceptional matrix-valued polynomials satisfying five- and seven-term
recurrence relations. We close by pointing out some directions for future
work, in particular the combination of matrix-valued objects and Dunkl
operators.

\section{Scalar ad-conditions}\label{sec:scalar}

In this section we establish the basic mechanism in the scalar case: a
sequence of polynomials which are eigenfunctions of a second-order
differential operator satisfies a (higher-order) recurrence relation if
and only if a certain operator identity, the \emph{ad-condition}, holds.

Throughout the paper $\partial$ denotes $\frac{d}{dx}$, and a second-order
differential operator with scalar coefficients is written as
\begin{equation}\label{eq:D-scalar}
    D = f_{2}(x)\,\partial^{2} + f_{1}(x)\,\partial + f_{0}(x),
    \qquad
    D\cdot p = f_{2}\,p'' + f_{1}\,p' + f_{0}\,p .
\end{equation}
Given an operator $\Theta$ (typically multiplication by a polynomial
$\Theta(x)$), the \emph{ad-ladder} of $\Theta$ with respect to $D$ is the
sequence of operators
\begin{equation}\label{eq:ladder}
    A_{0} = \Theta, \qquad
    A_{j} = \operatorname{ad}_{D}^{j}(\Theta) = [D,A_{j-1}] = [D,[D,\dots,[D,\Theta]\dots]],
    \quad j\in\mathbb{N} .
\end{equation}
The starting point is the following observation of M.~Reach \cite{R1}. We
state it in the generality in which it will be used later on: it never uses
that $D$ is a differential operator, nor that $\Theta$ is multiplication by
a function.

\begin{lemma}\label{lem:reach}
Let $\mathcal{A}$ be an associative unital $\mathbb{C}$-algebra acting on a vector space
$\mathcal{U}$, let $D,\Theta\in\mathcal{A}$, and let $A_{j}$ be as in \eqref{eq:ladder}. If
$\phi\in\mathcal{U}$ satisfies $D\cdot\phi=\lambda\phi$ with $\lambda\in\mathbb{C}$, then
for every $j\in\mathbb{N}_{0}$
\[
    (D-\lambda)^{j}\cdot(\Theta\cdot\phi)=A_{j}\cdot\phi .
\]
\end{lemma}

\begin{proof}
For $j=0$ there is nothing to prove. Assuming the statement for $j$, and
using that $\lambda$ is a scalar,
\[
    (D-\lambda)^{j+1}\cdot(\Theta\cdot\phi)=(D-\lambda)\cdot(A_{j}\cdot\phi)
    =D A_{j}\cdot\phi-A_{j}D\cdot\phi=[D,A_{j}]\cdot\phi=A_{j+1}\cdot\phi .\qedhere
\]
\end{proof}

\subsection{Three-term recurrence relations}\label{sec:three-term}

We fix $\Theta=x$, and let $\{p_{n}\}_{n\in\mathbb{N}_{0}}$ be a sequence of monic
polynomials with $\deg(p_{n})=n$. We say that $\{p_{n}\}$ satisfies a
three-term recurrence relation if
\[
    x\,p_{n}(x)=p_{n+1}(x)+b_{n}\,p_{n}(x)+c_{n}\,p_{n-1}(x),
    \qquad b_{n},c_{n}\in\mathbb{C} .
\]

\begin{theorem}\label{thm:three-term}
Let $D$ be a second-order differential operator \eqref{eq:D-scalar} such
that $D\cdot p_{n}=\lambda_{n}\,p_{n}$, $\lambda_{n}\in\mathbb{C}$, with
$\lambda_{n}\neq\lambda_{m}$ for all $n\neq m$. Then $\{p_{n}\}$ satisfies
a three-term recurrence relation if and only if for every $n\in\mathbb{N}_{0} $ one has
\begin{equation}\label{eq:cond}
    A_{3}\cdot p_{n}
    +(2\lambda_{n}-\lambda_{n-1}-\lambda_{n+1})\,A_{2}\cdot p_{n}
    +(\lambda_{n}-\lambda_{n+1})(\lambda_{n}-\lambda_{n-1})\,A_{1}\cdot p_{n}=0.
\end{equation}
\end{theorem}

\begin{proof}
First, we note that
\begin{equation}\label{eq:expansion3}
    \begin{split}
        & (D-\lambda_{n+1})(D-\lambda_{n})(D-\lambda_{n-1})\cdot(x\,p_{n}) \\
        & \quad = \bigl((D-\lambda_{n})+(\lambda_{n}-\lambda_{n+1})\bigr)(D-\lambda_{n})
            \bigl((D-\lambda_{n})+(\lambda_{n}-\lambda_{n-1})\bigr)\cdot(x\,p_{n}) \\
        & \quad = \bigl((D-\lambda_{n})^{3}+(2\lambda_{n}-\lambda_{n-1}-\lambda_{n+1})(D-\lambda_{n})^{2} \\
        & \qquad\qquad +(\lambda_{n}-\lambda_{n+1})(\lambda_{n}-\lambda_{n-1})(D-\lambda_{n})\bigr)\cdot(x\,p_{n}) \\
        & \quad = A_{3}\cdot p_{n}+(2\lambda_{n}-\lambda_{n-1}-\lambda_{n+1})\,A_{2}\cdot p_{n}
            +(\lambda_{n}-\lambda_{n+1})(\lambda_{n}-\lambda_{n-1})\,A_{1}\cdot p_{n},
    \end{split}
\end{equation}
where the last equality is Lemma~\ref{lem:reach}.

If $\{p_{n}\}$ satisfies a three-term recurrence relation,
$x\,p_{n}=p_{n+1}+b_{n}p_{n}+c_{n}p_{n-1}$, then clearly
\[
    (D-\lambda_{n+1})(D-\lambda_{n})(D-\lambda_{n-1})\cdot(x\,p_{n})=0 ,
\]
and \eqref{eq:cond} holds by \eqref{eq:expansion3}.

Conversely, assume that \eqref{eq:cond} holds. We write
\[
    x\,p_{n}(x)=p_{n+1}(x)+a_{n}p_{n}(x)+\cdots+a_{0}p_{0}(x),
    \qquad a_{0},\dots,a_{n}\in\mathbb{C} .
\]
From \eqref{eq:expansion3} we have that
\begin{equation*}
    \begin{split}
        0 & = (D-\lambda_{n+1})(D-\lambda_{n})(D-\lambda_{n-1})\cdot(x\,p_{n}) \\
          & = \sum_{j=0}^{n+1}(\lambda_{j}-\lambda_{n+1})(\lambda_{j}-\lambda_{n})(\lambda_{j}-\lambda_{n-1})\,a_{j}\,p_{j}(x) \\
          & = \sum_{j=0}^{n-2}(\lambda_{j}-\lambda_{n+1})(\lambda_{j}-\lambda_{n})(\lambda_{j}-\lambda_{n-1})\,a_{j}\,p_{j}(x).
    \end{split}
\end{equation*}
Since $\lambda_{j}\neq\lambda_{i}$ for all $j\neq i$, we obtain $a_{j}=0$
for $j=0,\dots,n-2$. Thus $\{p_{n}\}$ satisfies a three-term recurrence
relation.
\end{proof}

\begin{remark}\label{rem:even}
If $p_{2n}(x)$ is an even function and $p_{2n-1}(x)$ is an odd function
(as happens for the orthogonal polynomials with respect to an even weight
$w(x)$, for instance the Hermite and Gegenbauer polynomials), then
\eqref{eq:cond} can be replaced by the lower-order condition
\[
    A_{2}\cdot p_{n}
    +(2\lambda_{n}-\lambda_{n-1}-\lambda_{n+1})\,A_{1}\cdot p_{n}
    +(\lambda_{n}-\lambda_{n+1})(\lambda_{n}-\lambda_{n-1})\,A_{0}\cdot p_{n}=0 ,
\]
from which \eqref{eq:cond} follows by applying $\operatorname{ad}_{D}$.
\end{remark}

\begin{example}[Jacobi polynomials]\label{ex:jacobi}
Let $\{p_{n}\}$ be the sequence of monic Jacobi polynomials with parameters
$\alpha=\beta=2$. They are eigenfunctions of the operator
\[
    D=(1-x^{2})\,\partial^{2}-6x\,\partial,
    \qquad D\cdot p_{n}=\lambda_{n}\,p_{n},\quad \lambda_{n}=-(n^{2}+5n).
\]
Here $2\lambda_{n}-\lambda_{n-1}-\lambda_{n+1}=2$ and
$(\lambda_{n}-\lambda_{n+1})(\lambda_{n}-\lambda_{n-1})=-4(n+3)(n+2)=4(\lambda_{n}-6)$,
so that Theorem~\ref{thm:three-term} gives the ad-condition
\[
    \bigl(A_{3}+2A_{2}-4(n+3)(n+2)A_{1}\bigr)\cdot p_{n}=0 ,
    \qquad\text{that is,}\qquad
    \bigl(A_{3}+2A_{2}+4A_{1}(D-6)\bigr)\cdot p_{n}=0,
\]
for all $n\in\mathbb{N}_{0}$. The operator $A_{3}+2A_{2}+4A_{1}(D-6)$ is a
differential operator of order at most three which annihilates every
polynomial; since the kernel of a nonzero differential operator of order
$r$ has dimension at most $r$, it vanishes identically:
\[
    A_{3}+2A_{2}+4A_{1}(D-6)=0 .
\]
Replacing the $n$-dependent coefficient $-4(n+3)(n+2)$ by the operator
$4(D-6)$ is the device used by Reach in \cite[Case~II]{R1}.
\end{example}

\subsection{Higher-order recurrence relations and exceptional polynomials}\label{sec:higher}

We now allow for a finite set of \emph{exceptional degrees}, i.e.\ degrees
which do not occur in the sequence, and for a polynomial $\Theta$ of degree
$m\ge1$. Let $X$ be a finite subset of $\mathbb{N}_{0}$, the set of exceptional
degrees, and let $\{p_{n}\}$ be a sequence of monic polynomials with
$\deg(p_{n})=n$ for $n\in\mathbb{N}_{0}\setminus X$, and $p_{n}=0$ for $n\in X$. We
put
\[
    \mathcal{U}=\operatorname{span}_{\mathbb{C}}\{\,p_{n}(x)\,:\,n\in\mathbb{N}_{0}\setminus X\,\}.
\]
Among the sequences of this kind one has the so-called \emph{exceptional} orthogonal
polynomials \cite{GKM09,GKM10}, which arise naturally through Darboux transformations, see
\cite{GH96}. We say that $\{p_{n}\}$ satisfies a $(2m+1)$-term recurrence
relation with respect to a polynomial $\Theta(x)$ if
\[
    \Theta(x)\,p_{n}(x)=\sum_{i=-m}^{m}a_{n,i}\,p_{n+i}(x),
    \qquad a_{n,i}\in\mathbb{C},\quad n\in\mathbb{N}_{0}\setminus X .
\]
Exceptional orthogonal polynomials are known to satisfy recurrence relations
of this kind, see \cite{Dur15,GKKM16}.
Let $D$ be a second-order differential operator \eqref{eq:D-scalar}, and
assume
\[
    D\cdot p_{n}=\lambda_{n}\,p_{n},\qquad \lambda_{n}\in\mathbb{C},
    \quad n\in\mathbb{N}_{0}\setminus X,
\]
with $\lambda_{n}\neq\lambda_{m}$ if $n\neq m$. For an index $n+i$ with
$n+i<0$ or $n+i\in X$ the value $\lambda_{n+i}$ is not determined by the
eigenvalue equation; we assign to it an arbitrary value, subject only to
the condition that the $\lambda$'s remain pairwise distinct. In the
examples below $\lambda_{n}$ is given by a formula which makes sense for
all $n\in\mathbb{Z}$, and we use it.

\begin{theorem}\label{thm:BGZ-general}
In the situation above, $\{p_{n}\}$ satisfies a $(2m+1)$-term recurrence
relation with respect to $\Theta$ if and only if $\Theta\,\mathcal{U}\subset\mathcal{U}$ and
the ad-condition
\begin{equation}\label{eq:adcond-general}
    \Bigl(\sum_{j=1}^{2m+1}V_{j}(n)\,A_{j}\Bigr)\cdot p_{n}(x)=0,
    \qquad n\in\mathbb{N}_{0}\setminus X,
\end{equation}
holds, where $A_{j}$ is the ad-ladder \eqref{eq:ladder} of $\Theta$,
$V_{2m+1}(n)=1$, and
\begin{equation}\label{eq:Vj}
    V_{j}(n)=\sum_{-m\le i_{1}<i_{2}<\dots<i_{2m+1-j}\le m}\;
    \prod_{\ell=1}^{2m+1-j}\bigl(\lambda_{n}-\lambda_{n+i_{\ell}}\bigr),
    \qquad 1\le j\le 2m .
\end{equation}
\end{theorem}

\begin{proof}
Expanding the product and using Lemma~\ref{lem:reach},
\begin{align*}
    (D-\lambda_{n+m})\cdots(D-\lambda_{n-m})\cdot(\Theta\,p_{n})
    &=\left(\sum_{j=1}^{2m+1}V_{j}(n)\,(D-\lambda_{n})^{j}\right)\cdot(\Theta\,p_{n})
    \\
    &=\left(\sum_{j=1}^{2m+1}V_{j}(n)\,A_{j}\right)\cdot p_{n} .
\end{align*}
The proof then proceeds exactly as that of Theorem~\ref{thm:three-term}:
if the recurrence relation holds, the left-hand side vanishes; conversely,
if $\Theta\,\mathcal{U}\subset\mathcal{U}$ we may write $\Theta\,p_{n}=\sum_{j\le n+m}a_{j}p_{j}$, since
$\deg(\Theta\,p_{n})=n+m$,
and the vanishing of the left-hand side together with the distinctness of
the $\lambda$'s forces $a_{j}=0$ for $j<n-m$.
\end{proof}

\begin{example}[Exceptional Laguerre polynomials, five-term recurrence]\label{ex:laguerre-scalar}
Let $\alpha>1$, let $L_{n}^{(\alpha)}(x)$ denote the $n$-th monic Laguerre
polynomial of parameter $\alpha$, and consider the operator
\[
    \delta=-(\alpha+x)\,\partial^{2}+(\alpha+1+x)\,\partial .
\]
The polynomials
\[
    p_{0}^{(\alpha)}(x)=0,\qquad
    p_{n}^{(\alpha)}(x)=\frac{1}{n}\,\delta\cdot L_{n}^{(\alpha-2)}(x),
    \quad n\ge1,
\]
are monic of degree $n$, orthogonal with respect to the weight
\[
    w(x)=\frac{e^{-x}x^{\alpha}}{(x+\alpha)^{2}},\qquad x>0,
\]
(they are the $X_{1}$-Laguerre polynomials of \cite{GKM09}),
and they are eigenfunctions of the second-order operator
\begin{equation}\label{eq:D-laguerre}
    D=-x\,\partial^{2}
      -\frac{(\alpha-x)(x+\alpha+1)}{\alpha+x}\,\partial
      +\frac{2\alpha}{\alpha+x},
    \qquad D\cdot p_{n}^{(\alpha)}=n\,p_{n}^{(\alpha)} .
\end{equation}
Here $X=\{0\}$ and $\lambda_{n}=n$. For $\Theta(x)=(x+\alpha)^{2}$ the
sequence satisfies the five-term recurrence relation (cf.\ \cite{Dur15})
\begin{equation*}
    \begin{split}
        (x+\alpha)^{2}\,p_{n}^{(\alpha)}(x)
        =& p_{n+2}^{(\alpha)}(x)+4(\alpha+n)\,p_{n+1}^{(\alpha)}(x)
        \\
             &+\bigl(4\alpha^{2}+5(2n-1)\alpha+6n(n-1)\bigr)\,p_{n}^{(\alpha)}(x) \\
             &+4(n-1)(\alpha+n)(\alpha+n-2)\,p_{n-1}^{(\alpha)}(x) \\
             & +(n-1)(n-2)(\alpha+n)(\alpha+n-3)\,p_{n-2}^{(\alpha)}(x).
    \end{split}
\end{equation*}
Since the eigenvalue differences $\lambda_{n}-\lambda_{n+i}=-i$, $-2\le i\le2$,
do not depend on $n$, it is easy to see that the coefficients \eqref{eq:Vj} are given by $V_{5}=1$, $V_{4}=0$,
$V_{3}=-5$, $V_{2}=0$, $V_{1}=4$, and Theorem~\ref{thm:BGZ-general} with
$m=2$ gives
\[
    (A_{5}-5A_{3}+4A_{1})\cdot p_{n}^{(\alpha)}(x)=0,\qquad n\ge1 .
\]
As in Example~\ref{ex:jacobi}, the operator $A_{5}-5A_{3}+4A_{1}$ has
order at most five and annihilates the infinite-dimensional space $\mathcal{U}$;
hence
\[
    A_{5}-5A_{3}+4A_{1}=0 .
\]
\end{example}
\section{Ad-conditions in the Dunkl setting}\label{sec:dunkl}

All the ad-conditions obtained so far are identities between differential
operators, with coefficients that are scalars or polynomials in $n$. In this
section we show that the whole mechanism survives when $D$ is allowed to be a
\emph{Dunkl} operator, i.e.\ when the underlying algebra contains a reflection.
Two things happen, and both are new:

\begin{enumerate}
\item[(i)] Reach's original argument breaks down. In \cite{R1} the starting
point is that multiplication by $\Theta$ is an operator of order zero, so that
$\operatorname{ad}_{\Theta}^{3}(D)=0$ for $D$ of order two; this is what forces
$\lambda_{n}$ to be a polynomial of degree at most two in $n$. In the Dunkl algebra this
vanishing \emph{fails} (see Remark~\ref{rem:adtheta}), and accordingly the
eigenvalues $\lambda_{n}$ need not be polynomial in $n$.
\item[(ii)] Precisely because $\lambda_{n}$ is no longer polynomial in $n$, the
coefficients $V_{j}(n)$ of Theorem~\ref{thm:BGZ-general} cannot be scalars. They
are instead realized as \emph{elements of the algebra} acting diagonally on the
sequence $p_{n}$. In the rank-one Dunkl case the reflection $R$ does the job,
since $R\cdot p_{n}=(-1)^{n}p_{n}$, and one obtains ad-conditions of the
shape
\[
    A_{3}+c_{2}\,A_{2}R+c_{1}\,A_{1}=0 ,\qquad c_{1},c_{2}\in\mathbb{C},
\]
which are honest operator identities, with the reflection appearing as a
coefficient. This is the exact analogue of Reach's Case~II device of replacing
$n^{2}$ by the operator $D$ itself, see Example~\ref{ex:jacobi}.
\end{enumerate}

\subsection{The rank-one Dunkl algebra}

Let $R$ denote the reflection $R\cdot f(x)=f(-x)$, and for $\mu\in\mathbb{C}$ let
$T_{\mu}$ be the Dunkl operator \cite{Dun89}
\begin{equation}\label{eq:dunkl-op}
    T_{\mu}\cdot f(x)=f'(x)+\mu\,\frac{f(x)-f(-x)}{x}.
\end{equation}
We write $\mathcal{A}_{\mu}$ for the associative unital algebra generated by $T_{\mu}$,
the multiplication operator $x$, and $R$. Its defining relations are
\begin{equation}\label{eq:cherednik}
    [T_{\mu},x]=1+2\mu R,\qquad Rx=-xR,\qquad RT_{\mu}=-T_{\mu}R,\qquad R^{2}=1 .
\end{equation}
Note that $T_{\mu}$ lowers degrees by one, but with a parity-dependent factor:
on monomials
\[
    T_{\mu}\cdot x^{n}=\gamma_{n}\,x^{n-1},\qquad
    \gamma_{n}:=n+\mu\bigl(1-(-1)^{n}\bigr)
    =\begin{cases} n, & n \text{ even},\\[2pt] n+2\mu, & n \text{ odd}.\end{cases}
\]
The parity-dependent factor $\gamma_{n}$ is the source of everything that
follows.

\subsection{Spectral realizability}\label{sec:realizable}

Lemma~\ref{lem:reach} was stated for an arbitrary associative algebra
$\mathcal{A}$: it never uses that $D$ is a differential operator, nor that $\Theta$
is multiplication by a function. This is what allows us to take
$D,\Theta\in\mathcal{A}_{\mu}$. By contrast, the vanishing used by Reach is
\emph{not} algebra-independent.

\begin{remark}\label{rem:adtheta}
The vanishing $\operatorname{ad}_{\Theta}^{3}(D)=0$ used in
\cite[\S4]{R1} fails in $\mathcal{A}_{\mu}$. For $\Theta=x$ and
$D=D_{\mu}$ as in \eqref{eq:dunkl-hermite} below, a direct computation in
$\mathcal{A}_{\mu}$ using \eqref{eq:cherednik} gives
\[
    \operatorname{ad}_{x}^{3}(D_{\mu})=8\mu\,xR+16\mu\,x^{3}R,
\]
which is nonzero for $\mu\neq0$. Consequently the conclusion of \cite[\S5]{R1}
--- that $\lambda_{n}$ must be a polynomial of degree at most two in $n$ --- does
not apply here, and indeed the eigenvalues in \eqref{eq:lambda-hermite} are
\emph{not} polynomial in $n$. This is exactly the room that the Dunkl setting
creates.
\end{remark}

We now record the version of Theorem~\ref{thm:three-term} adapted to this
situation. The only change is that the coefficients $V_{j}(n)$ are required
to be the eigenvalues of operators $\mathcal{V}_{j}\in\mathcal{A}$, which are then placed to
the \emph{right} of $A_{j}$. In what follows $\mathcal{A}$ is an associative unital
algebra of linear operators on $\mathcal{U}=\mathbb{C}[x]$ (for instance $\mathcal{A}_{\mu}$).

\begin{definition}\label{def:realizable}
Let $\{p_{n}\}_{n\in\mathbb{N}_{0}}$ be a sequence in $\mathcal{U}$ with
$D\cdot p_{n}=\lambda_{n}p_{n}$, $D\in\mathcal{A}$. We say that a function
$V\colon\mathbb{N}_{0}\to\mathbb{C}$ is \emph{spectrally realizable} if there
exists $\mathcal{V}\in\mathcal{A}$ with
\[
    \mathcal{V}\cdot p_{n}=V(n)\,p_{n}\qquad\text{for all }n .
\]
\end{definition}

In the classical case of Section~\ref{sec:scalar} every polynomial in
$\lambda_{n}$ is spectrally realizable, by $\mathcal{V}=$ the corresponding polynomial in
$D$; this is why scalar coefficients suffice there, see Example~\ref{ex:jacobi}.
In the Dunkl case the algebra of realizable functions is strictly larger: it
contains $(-1)^{n}$, realized by $R$.

\begin{theorem}\label{thm:dunkl-adcond}
Let $\{p_{n}\}_{n\in\mathbb{N}_{0}}$ be a sequence of polynomials with
$\deg(p_{n})=n$, and let $D,\Theta\in\mathcal{A}$ with
\[
    D\cdot p_{n}(x)=\lambda_{n}p_{n}(x),\qquad
    \deg\bigl(\Theta\cdot p_{n}\bigr)\le n+1 ,
\]
and $\lambda_{j}\neq\lambda_{i}$ for all $j\neq i$. Assume that
\[
    V_{2}(n)=2\lambda_{n}-\lambda_{n-1}-\lambda_{n+1},\qquad
    V_{1}(n)=(\lambda_{n}-\lambda_{n+1})(\lambda_{n}-\lambda_{n-1})
\]
are spectrally realizable by $\mathcal{V}_{2},\mathcal{V}_{1}\in\mathcal{A}$. Then $\{p_{n}\}$ satisfies a
three-term recurrence relation with respect to $\Theta$, i.e.
\[
    \Theta\cdot p_{n}=a_{n}p_{n+1}+b_{n}p_{n}+c_{n}p_{n-1},
    \qquad a_{n},b_{n},c_{n}\in\mathbb{C},
\]
if and only if the \emph{ad-condition}
\begin{equation}\label{eq:dunkl-adcond}
    \bigl(A_{3}+A_{2}\mathcal{V}_{2}+A_{1}\mathcal{V}_{1}\bigr)\cdot p_{n}(x)=0
\end{equation}
holds for all $n$.
\end{theorem}

\begin{proof}
Exactly as in Theorem~\ref{thm:three-term}. By Lemma~\ref{lem:reach} and the
expansion
\[
    (D-\lambda_{n+1})(D-\lambda_{n})(D-\lambda_{n-1})
    =(D-\lambda_{n})^{3}+V_{2}(n)(D-\lambda_{n})^{2}+V_{1}(n)(D-\lambda_{n}),
\]
one gets
\[
    (D-\lambda_{n+1})(D-\lambda_{n})(D-\lambda_{n-1})\cdot\bigl(\Theta\cdot p_{n}\bigr)
    =\bigl(A_{3}+V_{2}(n)A_{2}+V_{1}(n)A_{1}\bigr)\cdot p_{n},
\]
and $V_{j}(n)A_{j}\cdot p_{n}=A_{j}\mathcal{V}_{j}\cdot p_{n}$ by spectral realizability.
If the three-term recurrence holds, the left-hand side vanishes. Conversely,
writing $\Theta\cdot p_{n}=\sum_{j\le n+1}a_{j}p_{j}$, the vanishing of
\eqref{eq:dunkl-adcond} gives
\[
    0=\sum_{j=0}^{n-2}(\lambda_{j}-\lambda_{n+1})(\lambda_{j}-\lambda_{n})
      (\lambda_{j}-\lambda_{n-1})\,a_{j}\,p_{j}(x),
\]
and distinctness of the $\lambda_{j}$ forces $a_{j}=0$ for $j\le n-2$.
\end{proof}

\begin{remark}
The order of the factors in \eqref{eq:dunkl-adcond} matters: $\mathcal{V}_{j}$ must be
written on the right of $A_{j}$, since in general $[\,A_{j},\mathcal{V}_{j}\,]\neq0$.
When $\mathcal{V}_{j}$ is a scalar one recovers \eqref{eq:cond} verbatim.
\end{remark}

\subsection{Generalized Hermite polynomials}\label{sec:gen-hermite}

\begin{example}[Generalized Hermite polynomials]\label{ex:gen-hermite}
Let $\mu>-\tfrac12$ and consider the weight
\[
    w_{\mu}(x)=|x|^{2\mu}e^{-x^{2}},\qquad x\in\mathbb{R},
\]
with associated sequence $\{H^{\mu}_{n}\}_{n\ge0}$ of monic orthogonal polynomials, the \emph{generalized Hermite polynomials}, first considered by
G.~Szeg\H{o} \cite[p.~380, Problem~25]{Sze75}, see also \cite[Chapter~V, \S2]{Chi78}. They satisfy
\begin{equation}\label{eq:gh-rec}
    x\,H^{\mu}_{n}(x)=H^{\mu}_{n+1}(x)+\tfrac{1}{2}\gamma_{n}\,H^{\mu}_{n-1}(x),
    \qquad \gamma_{n}=n+\mu\bigl(1-(-1)^{n}\bigr),
\end{equation}
and, since $w_{\mu}$ is even, $H^{\mu}_{n}(-x)=(-1)^{n}H^{\mu}_{n}(x)$, that is
\begin{equation}\label{eq:parity}
    R\cdot H^{\mu}_{n}=(-1)^{n}H^{\mu}_{n}.
\end{equation}
Consider the \emph{Dunkl harmonic oscillator} \cite{Ros94}
\begin{equation}\label{eq:dunkl-hermite}
    D_{\mu}=T_{\mu}^{2}-2x\,T_{\mu}\ \in\ \mathcal{A}_{\mu}.
\end{equation}
Then $D_{\mu}$ preserves degrees and, see \cite{Ros94},
\begin{equation}\label{eq:lambda-hermite}
    D_{\mu}\cdot H^{\mu}_{n}(x)=\lambda_{n}H^{\mu}_{n}(x),
    \qquad \lambda_{n}=-2\gamma_{n}
    =\begin{cases} -2n, & n\text{ even},\\[2pt] -2n-4\mu, & n\text{ odd}.\end{cases}
\end{equation}
\end{example}

\begin{remark}\label{rem:distinct}
The eigenvalues $\lambda_{n}$ are pairwise distinct if and only if
$\mu\notin\tfrac12+\mathbb{Z}$. Indeed $\lambda_{n}=\lambda_{m}$ with $n$ even
and $m$ odd forces $\mu=(n-m)/2$ with $n-m$ odd. Together with the
integrability condition this gives the admissible range
\[
    \mu\in\Bigl(-\tfrac12,\infty\Bigr)\setminus
    \Bigl\{\tfrac12,\tfrac32,\tfrac52,\dots\Bigr\}.
\]
Only the converse half of Theorem~\ref{thm:dunkl-adcond} needs this hypothesis;
the ad-condition itself holds for every $\mu$.
\end{remark}

We now compute the coefficients. From \eqref{eq:lambda-hermite},
\[
    \lambda_{n}-\lambda_{n+1}=2\bigl(1+2\mu(-1)^{n}\bigr),\qquad
    \lambda_{n}-\lambda_{n-1}=-2\bigl(1-2\mu(-1)^{n}\bigr),
\]
so that
\begin{equation}\label{eq:V-hermite}
    V_{2}(n)=2\lambda_{n}-\lambda_{n-1}-\lambda_{n+1}=8\mu(-1)^{n},
    \qquad
    V_{1}(n)=(\lambda_{n}-\lambda_{n+1})(\lambda_{n}-\lambda_{n-1})=-4\bigl(1-4\mu^{2}\bigr).
\end{equation}
The second coefficient is a constant, but the first one is \emph{not} a
polynomial in $n$ --- it is not even bounded away from oscillating. By
\eqref{eq:parity} it is nonetheless spectrally realizable:
\[
    \mathcal{V}_{2}=8\mu R,\qquad \mathcal{V}_{1}=-4(1-4\mu^{2}) .
\]
Note also that $[D_{\mu},R]=0$, so $R$ commutes with the whole ad-ladder
construction.

\begin{theorem}\label{thm:hermite-adcond}
Let $\Theta=x$ and $A_{j}=\operatorname{ad}_{D_{\mu}}^{j}(x)$. Then the following
identity holds in $\mathcal{A}_{\mu}$, for every $\mu\in\mathbb{C}$:
\begin{equation}\label{eq:hermite-adcond}
    A_{3}+8\mu\,A_{2}R-4\bigl(1-4\mu^{2}\bigr)A_{1}=0 .
\end{equation}
If moreover $\mu\notin\tfrac12+\mathbb{Z}$, then \eqref{eq:hermite-adcond} is
equivalent to the three-term recurrence relation \eqref{eq:gh-rec}.
\end{theorem}

\begin{proof}
The second assertion is Theorem~\ref{thm:dunkl-adcond} together with
\eqref{eq:V-hermite} and Remark~\ref{rem:distinct}. For the first, we compute the
ad-ladder explicitly. Using \eqref{eq:cherednik} and $RT_{\mu}=-T_{\mu}R$,
\[
    [T_{\mu}^{2},x]=T_{\mu}(1+2\mu R)+(1+2\mu R)T_{\mu}
    =2T_{\mu}+2\mu\,(T_{\mu}R+RT_{\mu})=2T_{\mu},
\]
\[
    [-2xT_{\mu},x]=-2x[T_{\mu},x]=-2x-4\mu\,xR,
\]
so that
\begin{equation}\label{eq:A1}
    A_{1}=2T_{\mu}-2x-4\mu\,xR .
\end{equation}
Next, $[D_{\mu},R]=0$ (both $T_{\mu}^{2}$ and $xT_{\mu}$ commute with $R$) and
\begin{equation}\label{eq:DT}
    [D_{\mu},T_{\mu}]=-2[xT_{\mu},T_{\mu}]=2(1+2\mu R)T_{\mu}=2T_{\mu}+4\mu RT_{\mu},
\end{equation}
whence $[D_{\mu},xR]=A_{1}R$ and
\begin{equation}\label{eq:A2}
    A_{2}=2[D_{\mu},T_{\mu}]-2A_{1}-4\mu A_{1}R
        =4x+16\mu\,xR+16\mu^{2}x-16\mu\,T_{\mu}R .
\end{equation}
Applying $\operatorname{ad}_{D_{\mu}}$ once more, and using
$[D_{\mu},T_{\mu}R]=2T_{\mu}R-4\mu T_{\mu}$ from \eqref{eq:DT},
\begin{equation}\label{eq:A3}
    A_{3}=\bigl(8+96\mu^{2}\bigr)\bigl(T_{\mu}-x\bigr)
          -\bigl(48\mu+64\mu^{3}\bigr)xR .
\end{equation}
Finally, from \eqref{eq:A2} and $R^{2}=1$,
\[
    8\mu A_{2}R=32\mu\,xR+128\mu^{2}x+128\mu^{3}xR-128\mu^{2}T_{\mu},
\]
and from \eqref{eq:A1},
\[
    -4(1-4\mu^{2})A_{1}=-\bigl(8-32\mu^{2}\bigr)T_{\mu}+\bigl(8-32\mu^{2}\bigr)x
    +\bigl(16\mu-64\mu^{3}\bigr)xR .
\]
Adding the three expressions, the coefficients of $T_{\mu}$, of $x$ and of $xR$
are respectively
\[
    (8+96\mu^{2})-128\mu^{2}-(8-32\mu^{2})=0,\qquad
    -(8+96\mu^{2})+128\mu^{2}+(8-32\mu^{2})=0,
\]
\[
    -(48\mu+64\mu^{3})+32\mu+128\mu^{3}+16\mu-64\mu^{3}=0,
\]
while the coefficients of $T_{\mu}R$ cancel already in $A_{3}$. This proves
\eqref{eq:hermite-adcond}.
\end{proof}

In fact the ad-condition for this example descends one step further, which is
the Dunkl analogue of Remark~\ref{rem:even} for even weights.

\begin{proposition}\label{prop:minimal}
With the notation of Theorem~\ref{thm:hermite-adcond}, the identity
\begin{equation}\label{eq:hermite-minimal}
    A_{2}+8\mu\,A_{1}R-4\bigl(1-4\mu^{2}\bigr)A_{0}=0
\end{equation}
holds in $\mathcal{A}_{\mu}$, where $A_{0}=x$. Since $[D_{\mu},R]=0$, applying
$\operatorname{ad}_{D_{\mu}}$ to \eqref{eq:hermite-minimal} returns
\eqref{eq:hermite-adcond}. Thus \eqref{eq:hermite-minimal} is the minimal
ad-condition satisfied by the generalized Hermite polynomials, and
\eqref{eq:hermite-adcond} is a consequence of it.
\end{proposition}

\begin{proof}
From \eqref{eq:A1}, $8\mu A_{1}R=16\mu T_{\mu}R-16\mu\,xR-32\mu^{2}x$. Adding
\eqref{eq:A2} and $-4(1-4\mu^{2})x=-4x+16\mu^{2}x$, the coefficients of
$T_{\mu}R$, $xR$ and $x$ are
\[
    -16\mu+16\mu=0,\qquad 16\mu-16\mu=0,\qquad
    4+16\mu^{2}-32\mu^{2}-4+16\mu^{2}=0 .
\]
The last assertion follows since $\operatorname{ad}_{D_{\mu}}(A_{j}R)=A_{j+1}R$.
\end{proof}

\begin{remark}
For $\mu=0$ one has $T_{0}=\partial$, $D_{0}=\partial^{2}-2x\partial$,
and \eqref{eq:hermite-adcond} collapses to
\[
    A_{3}-4A_{1}=0,\qquad A_{1}=2\partial-2x,\quad A_{2}=4x,\quad A_{3}=4A_{1},
\]
which is Reach's Case~I condition $A_{3}-A_{1}=0$ after the rescaling
$D\mapsto-\tfrac12 D$ of \cite[(S1)]{R1}. For $\mu\neq0$ the middle term
$A_{2}$ is genuinely present, and \emph{no} choice of scalar coefficients can
produce it: by \eqref{eq:V-hermite} its coefficient alternates in sign with $n$.
The reflection is therefore not a convenience but a necessity. Operators in
which the reflection enters as a coefficient appear also in the bispectral
treatment of the Bannai--Ito polynomials by means of Dunkl shift operators,
see \cite{VZ11,TVZ12}.
\end{remark}

\begin{remark}
The two degenerations visible in \eqref{eq:hermite-adcond} coincide: the
coefficient $-4(1-4\mu^{2})$ of $A_{1}$ vanishes exactly at $\mu=\pm\tfrac12$,
and $\mu=\tfrac12$ is precisely the first value excluded in
Remark~\ref{rem:distinct}, where $\lambda_{2n+1}=\lambda_{2n+2}$. At $\mu=\tfrac12$
the ad-condition reduces to $A_{3}+4A_{2}R=0$, and the converse of
Theorem~\ref{thm:dunkl-adcond} is no longer available. The ad-condition
``sees'' the collision of the spectrum.
\end{remark}

\subsection{A five-term example: $\Theta=x^{2}$}\label{sec:x2}

\begin{example}\label{ex:x2}
We keep the notation of Example~\ref{ex:gen-hermite}. Taking $\Theta=x^{2}$
places us in the setting of Theorem~\ref{thm:BGZ-general} with $m=2$.
Iterating \eqref{eq:gh-rec},
\begin{equation}\label{eq:gh-five}
    x^{2}H^{\mu}_{n}(x)=H^{\mu}_{n+2}(x)
    +\tfrac{1}{2}\bigl(\gamma_{n+1}+\gamma_{n}\bigr)H^{\mu}_{n}(x)
    +\tfrac{1}{4}\gamma_{n}\gamma_{n-1}\,H^{\mu}_{n-2}(x),
\end{equation}
a five-term recurrence relation whose coefficients at $n\pm1$ vanish. This
sparsity is inherited by the ad-condition: instead of the generic product over
the five offsets $\{-2,-1,0,1,2\}$ it suffices to use the three offsets
$\{-2,0,2\}$. Since $\gamma_{n\pm2}=\gamma_{n}\pm2$, one has
\[
    \lambda_{n}-\lambda_{n\pm2}=\pm4,
\]
\emph{independently of $n$ and of $\mu$}, so
\[
    (D_{\mu}-\lambda_{n+2})(D_{\mu}-\lambda_{n})(D_{\mu}-\lambda_{n-2})
    =(D_{\mu}-\lambda_{n})^{3}-16\,(D_{\mu}-\lambda_{n}) .
\]
\end{example}

\begin{proposition}\label{prop:x2}
Let $C_{j}=\operatorname{ad}_{D_{\mu}}^{j}(x^{2})$. Then, for every
$\mu\in\mathbb{C}$,
\begin{equation}\label{eq:x2-adcond}
    C_{3}-16\,C_{1}=0 .
\end{equation}
Explicitly,
\[
    C_{1}=4xT_{\mu}+2+4\mu R-4x^{2},\qquad
    C_{2}=8T_{\mu}^{2}-16xT_{\mu}-8-16\mu R+16x^{2},\qquad C_{3}=16C_{1}.
\]
\end{proposition}

\begin{proof}
From \eqref{eq:cherednik}, $[T_{\mu},x^{2}]=[T_{\mu},x]x+x[T_{\mu},x]
=(1+2\mu R)x+x(1+2\mu R)=2x$, the reflection dropping out because $Rx=-xR$.
Hence
\[
    [T_{\mu}^{2},x^{2}]=2T_{\mu}x+2xT_{\mu}=4xT_{\mu}+2+4\mu R,\qquad
    [-2xT_{\mu},x^{2}]=-4x^{2},
\]
which gives $C_{1}$. Next, using \eqref{eq:A1} and \eqref{eq:DT},
\[
    [D_{\mu},xT_{\mu}]=A_{1}T_{\mu}+x[D_{\mu},T_{\mu}]
    =\bigl(2T_{\mu}-2x-4\mu xR\bigr)T_{\mu}+x\bigl(2T_{\mu}+4\mu RT_{\mu}\bigr)
    =2T_{\mu}^{2},
\]
so $C_{2}=4[D_{\mu},xT_{\mu}]-4C_{1}=8T_{\mu}^{2}-4C_{1}$. Finally
$[D_{\mu},T_{\mu}^{2}]=4T_{\mu}^{2}$ by \eqref{eq:DT} and $RT_{\mu}=-T_{\mu}R$,
so
\[
    C_{3}=8\,[D_{\mu},T_{\mu}^{2}]-4C_{2}=32T_{\mu}^{2}-4\bigl(8T_{\mu}^{2}-4C_{1}\bigr)=16C_{1}.\qedhere
\]
\end{proof}

\begin{remark}
The ad-condition \eqref{eq:x2-adcond} is free of both $R$ and $\mu$, in sharp
contrast with \eqref{eq:hermite-adcond}. The reason is structural: $\Theta=x^{2}$
is \emph{even}, so it preserves each of the two parity subspaces on which $R$ is
scalar, and the parity-dependent part of the spectrum becomes invisible. In the
language of Theorem~\ref{thm:dunkl-adcond}, all the relevant $V_{j}$ are already
scalars and no spectral realization is needed. The whole $\mu$-family
$\{D_{\mu}\}_{\mu}$ satisfies one and the same ad-condition with respect to
$x^{2}$.  
\end{remark}
 
One could also consider the generalized Gegenbauer polynomials. We refrain from doing so here
because of page limitations.

\section{Matrix-valued ad-conditions}\label{sec:matrix}

In this section we consider sequences of $N\times N$ matrix-valued
polynomials $\mathcal{P}_{n}(x)$ which are eigenfunctions of a matrix-valued
differential operator and satisfy a higher-order recurrence relation; in
all our examples $N=2$. We do not attempt a general theory here. After
fixing the conventions, we exhibit, for three families of exceptional
matrix-valued polynomials with five- and seven-term recurrence relations,
built from the scalar examples of Section~\ref{sec:scalar} by means of the
constructions of \cite{B26,BDZ25}, operator identities of exactly the same
shape as the scalar ad-conditions of Section~\ref{sec:scalar}. Each of
these identities is established by a direct computation.

\subsection{Conventions}\label{sec:matrix-conventions}

We follow the usual convention in the theory of matrix-valued orthogonal
polynomials, see \cite{GPT03,B26}: differential operators act on the
\emph{right},
\begin{equation}\label{eq:D-matrix}
    \mathcal{D}=\partial^{2}F_{2}(x)+\partial F_{1}(x)+F_{0}(x),
    \qquad
    \mathcal{P}\cdot\mathcal{D}=\mathcal{P}''F_{2}+\mathcal{P}'F_{1}+\mathcal{P} F_{0},
\end{equation}
so that the eigenvalue equation takes the form
\begin{equation}\label{eq:eig-matrix}
    \mathcal{P}_{n}\cdot\mathcal{D}=\Lambda_{n}\,\mathcal{P}_{n},
\end{equation}
with a constant matrix $\Lambda_{n}$ on the left, and a $(2m+1)$-term
recurrence relation takes the form
\begin{equation}\label{eq:rec-matrix}
    \mathcal{P}_{n}(x)\,\Theta(x)=\sum_{i=-m}^{m}B_{i}(n)\,\mathcal{P}_{n+i}(x),
\end{equation}
with a matrix-valued polynomial $\Theta(x)$ on the right and constant
matrices $B_{i}(n)$ on the left. For scalar operators the two conventions
give the same action, and the scalar operators of Section~\ref{sec:scalar}
reappear below acting on the right. Products of right-acting operators are
read from left to right: $\mathcal{P}\cdot(XY)=(\mathcal{P}\cdot X)\cdot Y$. With this
convention the ad-ladder of $\Theta$ with respect to $\mathcal{D}$ is still
defined by \eqref{eq:ladder}, $A_{0}=\Theta$, $A_{j}=[\mathcal{D},A_{j-1}]$. Note
that reading products from left to right reverses the order of
composition with respect to Section~\ref{sec:scalar}, so that for a scalar
$D$ and a scalar $\Theta$ the ladder of this section is $(-1)^{j}$ times
the ladder of Section~\ref{sec:scalar}. All the identities below involve
only odd, or only even, values of $j$, and are therefore insensitive to
this sign.

In each of the examples below the matrices $\Lambda_{n}$ are diagonal and
the differences $\Lambda_{n}-\Lambda_{n+i}$ are scalar and independent of
$n$, exactly as in the corresponding scalar examples of
Section~\ref{sec:scalar}. This suggests looking for an operator identity
$\sum_{j}V_{j}A_{j}=0$ with the same scalar coefficients $V_{j}$ as in
\eqref{eq:Vj}. In each case such an identity does hold, and it is the
ad-condition that we exhibit; in particular the operator
$\sum_{j}V_{j}A_{j}$ annihilates every $\mathcal{P}_{n}$.

\subsection{Matrix-valued examples}\label{sec:matrix-examples}

\begin{example}[Exceptional Laguerre polynomials, five-term recurrence]\label{ex:laguerre-5}
Let $\alpha>1$ and let $p_{n}^{(\alpha)}$, $w$ and $D$ be as in
Example~\ref{ex:laguerre-scalar}. Adapting \cite[Theorems 3.1 and 4.1]{B26},
we build, for $b\neq0$, the following sequence of matrix-valued exceptional
Laguerre polynomials with exceptional degree $0$,
\[
    \mathcal{P}_{n}(x)=\begin{pmatrix}
        p_{n}^{(\alpha)}(x) & b\,\bigl(p_{n+1}^{(\alpha)}(x)-x\,p_{n}^{(\alpha)}(x)\bigr) \\[4pt]
        -b\,\kappa_{n}\,p_{n-1}^{(\alpha)}(x) & b^{2}\,\kappa_{n}\,x\,p_{n-1}^{(\alpha)}(x)+p_{n}^{(\alpha)}(x)
    \end{pmatrix},
    \qquad n\ge1 ,
\]
with
\[
    \kappa_{n}:=\frac{\|p_{n}^{(\alpha)}\|^{2}}{\|p_{n-1}^{(\alpha)}\|^{2}}
    =\frac{(n+\alpha-2)(n-1)(n+\alpha)}{n+\alpha-1},
\]
where $\|p\|^{2}=\int_{0}^{\infty}p(x)^{2}w(x)\,dx$.
It is easy to verify that they are orthogonal with respect to the weight matrix
\[
    \mathcal{W}(x)=\frac{e^{-x}x^{\alpha}}{(x+\alpha)^{2}}
    \begin{pmatrix} 1+b^{2}x^{2} & bx \\ bx & 1\end{pmatrix},\qquad x>0 ,
\]
and are eigenfunctions of the second-order operator
\begin{equation*}
    \begin{split}
        \mathcal{D} = {} & \partial^{2}\begin{pmatrix}-x & 0 \\ 0 & -x\end{pmatrix}
        +\partial\begin{pmatrix}
            -\dfrac{(\alpha-x)(x+\alpha+1)}{\alpha+x} & -2bx \\[8pt]
            0 & -\dfrac{(\alpha-x)(x+\alpha+1)}{\alpha+x}
        \end{pmatrix} \\
        & +\begin{pmatrix}
            \dfrac{2\alpha}{\alpha+x} & -\dfrac{b\,(\alpha^{2}+\alpha x+\alpha-x)}{\alpha+x} \\[8pt]
            0 & \dfrac{\alpha-x}{\alpha+x}
        \end{pmatrix},
    \end{split}
\end{equation*}
with
\[
    \mathcal{P}_{n}\cdot\mathcal{D}=\begin{pmatrix} n & 0 \\ 0 & n-1\end{pmatrix}\mathcal{P}_{n},\qquad n\ge1 .
\]
Moreover, $\mathcal{P}_{n}$ satisfies a five-term recurrence relation
\eqref{eq:rec-matrix} with respect to the scalar polynomial
$\Theta(x)=x^{2}+2\alpha x$. The eigenvalue differences are
$\Lambda_{n}-\Lambda_{n+i}=-i\,I$, $-2\le i\le2$, exactly as in
Example~\ref{ex:laguerre-scalar}, and the coefficients \eqref{eq:Vj} are
again $V_{5}=1$, $V_{3}=-5$, $V_{1}=4$, with $V_{4}=V_{2}=0$. A direct
computation shows that the operator identity
\[
    A_{5}-5A_{3}+4A_{1}=0
\]
holds for the ad-ladder of $\Theta$ with respect to $\mathcal{D}$, and
consequently also $A_{6}-5A_{4}+4A_{2}=0$. In particular
$A_{5}-5A_{3}+4A_{1}$ annihilates every $\mathcal{P}_{n}$. Note that
$\Theta(x)$ differs from the polynomial $(x+\alpha)^{2}$ of
Example~\ref{ex:laguerre-scalar} by a constant, so that the two ad-ladders
agree for $j\ge1$, and the scalar operator $D$ satisfies the same
ad-condition.
\end{example}

\begin{example}[Exceptional Laguerre polynomials, seven-term recurrence]\label{ex:laguerre-7}
Let $\alpha>2$, and consider the first-order operator
\[
    \widetilde{\delta}=x\,L_{2}^{(-\alpha)}(x)\,\partial
    +\bigl(\alpha\,L_{2}^{(-\alpha)}(x)-2x\,L_{1}^{(-\alpha+1)}(x)\bigr).
\]
For $n\ge0$ the polynomials
\[
    p_{n+2}^{(\alpha)}(x)=\frac{1}{\alpha+n-2}\,\widetilde{\delta}\cdot L_{n}^{(\alpha)}(x)
\]
are monic of degree $n+2$, orthogonal with respect to the weight
\[
    w(x)=\frac{e^{-x}x^{\alpha-1}}{\bigl(L_{2}^{(-\alpha)}(x)\bigr)^{2}},\qquad x>0,
\]
and eigenfunctions of
\[
    D=-x\,\partial^{2}
      +\Bigl((x-\alpha)+\frac{4x\,L_{1}^{(-\alpha+1)}(x)}{L_{2}^{(-\alpha)}(x)}\Bigr)\partial
      +\alpha+1-\frac{4x\,L_{1}^{(-\alpha+1)}(x)}{L_{2}^{(-\alpha)}(x)},
\]
with
\[
    D\cdot p_{n+2}^{(\alpha)}=(\alpha+n-1)\,p_{n+2}^{(\alpha)},\qquad n\ge0 .
\]
Here the exceptional degrees are $\{0,1\}$. The sequence satisfies a
seven-term recurrence relation with respect to the polynomial
\[
    \Theta(x)=x^{3}+(3\alpha-6)x^{2}+(3\alpha^{2}-9\alpha+6)x ,
\]
and since $\lambda_{n}-\lambda_{n+i}=-i$, $-3\le i\le3$, the coefficients
\eqref{eq:Vj} are $V_{7}=1$, $V_{5}=-14$, $V_{3}=49$, $V_{1}=-36$, the even
ones vanishing. Theorem~\ref{thm:BGZ-general} with $m=3$ and the argument
of Example~\ref{ex:laguerre-scalar} give the scalar ad-condition
\[
    A_{7}-14A_{5}+49A_{3}-36A_{1}=0
\]
for the ad-ladder of $\Theta$ with respect to $D$.

As in the previous example, adapting \cite[Theorems 3.1 and 4.1]{B26} we
build, for $b\neq0$, the following sequence of matrix-valued exceptional
Laguerre polynomials, whose exceptional degrees are also $\{0,1\}$,
\[
    \mathcal{P}_{n+2}(x)=\begin{pmatrix}
        p_{n+2}^{(\alpha)}(x) & b\,\bigl(p_{n+3}^{(\alpha)}(x)-x\,p_{n+2}^{(\alpha)}(x)\bigr) \\[4pt]
        -b\,\rho_{n}\,p_{n+1}^{(\alpha)}(x) & b^{2}\,\rho_{n}\,x\,p_{n+1}^{(\alpha)}(x)+p_{n+2}^{(\alpha)}(x)
    \end{pmatrix},\qquad n\ge0,
\]
where
\[
    \rho_{n}:=\frac{\|p_{n+2}^{(\alpha)}\|^{2}}{\|p_{n+1}^{(\alpha)}\|^{2}}
    =\frac{n(n+\alpha)(n+\alpha-3)}{n+\alpha-2} ,
\]
the norms being taken with respect to the weight $w$ above.
It is easy to verify that they are orthogonal with respect to the weight matrix
\[
    \mathcal{W}(x)=\frac{e^{-x}x^{\alpha-1}}{\bigl(L_{2}^{(-\alpha)}(x)\bigr)^{2}}
    \begin{pmatrix} 1+b^{2}x^{2} & bx \\ bx & 1\end{pmatrix},\qquad x>0,
\]
and eigenfunctions of the second-order operator
\begin{equation*}
    \begin{split}
        \mathcal{D} = {} & \partial^{2}\begin{pmatrix}-x & 0 \\ 0 & -x\end{pmatrix}
        +\partial\begin{pmatrix}
            (x-\alpha)+\dfrac{4x\,L_{1}^{(-\alpha+1)}(x)}{L_{2}^{(-\alpha)}(x)} & -2bx \\[8pt]
            0 & (x-\alpha)+\dfrac{4x\,L_{1}^{(-\alpha+1)}(x)}{L_{2}^{(-\alpha)}(x)}
        \end{pmatrix} \\
        & +\begin{pmatrix}
            \alpha+1-\dfrac{4x\,L_{1}^{(-\alpha+1)}(x)}{L_{2}^{(-\alpha)}(x)}
            & -b\Bigl(\alpha-\dfrac{4x\,L_{1}^{(-\alpha+1)}(x)}{L_{2}^{(-\alpha)}(x)}\Bigr) \\[8pt]
            0 & \alpha-\dfrac{4x\,L_{1}^{(-\alpha+1)}(x)}{L_{2}^{(-\alpha)}(x)}
        \end{pmatrix},
    \end{split}
\end{equation*}
with
\[
    \mathcal{P}_{n+2}\cdot\mathcal{D}=\begin{pmatrix}\alpha+n-1 & 0 \\ 0 & \alpha+n-2\end{pmatrix}\mathcal{P}_{n+2},
    \qquad n\ge0 .
\]
The sequence $\mathcal{P}_{n}$ satisfies a seven-term recurrence relation
\eqref{eq:rec-matrix} with respect to the same scalar polynomial
$\Theta(x)$, and the eigenvalue differences are again
$\Lambda_{n}-\Lambda_{n+i}=-i\,I$, $-3\le i\le3$. A direct computation
shows that the operator identity
\[
    A_{7}-14A_{5}+49A_{3}-36A_{1}=0
\]
holds for the ad-ladder of $\Theta$ with respect to $\mathcal{D}$, so that the
matrix-valued operator $\mathcal{D}$ satisfies, with respect to $\Theta$, the same
ad-condition as the scalar operator $D$; in particular this operator
annihilates every $\mathcal{P}_{n}$.
\end{example}

\begin{example}[Exceptional Hermite polynomials, five-term recurrence]\label{ex:hermite-matrix}
Using the results of \cite{BDZ25} we construct a family of matrix-valued
exceptional Hermite polynomials. Let $H_{n}(x)$ be the monic Hermite
polynomials, let $a\in\mathbb{R}$ with $a^{2}>2$, and consider the
matrix-valued polynomials
\[
    P_{n}(x)=\begin{pmatrix}
        H_{n}(x) & -\dfrac{a}{2}\,n\,H_{n-1}(x) \\[6pt]
        -\dfrac{a}{2}\,n\,H_{n-1}(x) & (1-a^{2})H_{n}(x)+\dfrac{a^{2}}{2}\,n\,x\,H_{n-1}(x)
    \end{pmatrix}
\]
and the first-order operator
\[
    \Delta=\partial\begin{pmatrix} x & -\dfrac{a}{2} \\[6pt] \dfrac{a}{a^{2}-2} & x\end{pmatrix}-I .
\]
We define the sequence of exceptional matrix-valued polynomials
\[
    \mathcal{P}_{n}(x)=P_{n}\cdot\Delta ,
\]
which have degree $n$ for $n\neq1$ and $\mathcal P_1=0$. They are
eigenfunctions of the operator $\mathcal{D}=\partial^{2}+\partial F_{1}(x)+F_{0}(x)$,
with
\[
    F_{1}(x)=\frac{1}{q(x)}\begin{pmatrix}
        -2x\,(q(x)-4) & 2a(a^{2}-2)(2x^{2}-1) \\[6pt]
        -\dfrac{8a}{a^{2}-2} & -2x\,(q(x)-4+4a^{2})
    \end{pmatrix},
    \qquad
    F_{0}(x)=\begin{pmatrix} -2 & 0 \\ 0 & 0\end{pmatrix},
\]
where $q(x)=2a^{2}x^{2}+a^{2}-4x^{2}$, and
\[
    \mathcal{P}_{n}\cdot\mathcal{D}=\begin{pmatrix} -2(n+1) & 0 \\ 0 & -2n\end{pmatrix}\mathcal{P}_{n} .
\]
The polynomials $\mathcal{P}_{n}$ are orthogonal with respect to the weight matrix

\[
    \mathcal{W}(x)=\frac{e^{-x^{2}}}{q(x)^{2}}
    \begin{pmatrix}
        a^{2}q(x)^{2}-4q(x)-4a^{2}(a^{2}-2) & -4a x\bigl(q(x)+2(a^{2}-2)\bigr) \\[6pt]
        -4a\,x\,\bigl(q(x)+2(a^{2}-2)\bigr) & \dfrac{8\bigl(q(x)+a^{2}(a^{2}-2)\bigr)}{a^{2}-2}
    \end{pmatrix}
\]
 on the real line.
The sequence $\mathcal{P}_{n}$ also satisfies a five-term recurrence relation
\eqref{eq:rec-matrix} with respect to the matrix-valued polynomial
\[
    \Theta(x)=\begin{pmatrix}
        -\dfrac{a^{3}(a^{2}-2)}{2(a^{2}-1)}\,x & -\dfrac{(a^{2}-2)^{2}(2a^{2}x^{2}-a^{2}+2)}{8(a^{2}-1)} \\[8pt]
        1 & 0
    \end{pmatrix}.
\]
The eigenvalue differences are $\Lambda_{n}-\Lambda_{n+i}=2i\,I$, $-2\le i\le2$, so
that the coefficients \eqref{eq:Vj} are $V_{5}=1$, $V_{3}=-20$, $V_{1}=64$,
the even ones vanishing. A direct computation shows that the operator
identity
\[
    A_{5}-20A_{3}+64A_{1}=0
\]
holds for the ad-ladder of $\Theta$ with respect to $\mathcal{D}$, so that this
operator annihilates every $\mathcal{P}_{n}$, and that it descends one step
further, as in Remark~\ref{rem:even}, to
\[
    A_{4}-20A_{2}+64A_{0}=0 .
\]
\end{example}
\section{Concluding remarks}\label{sec:conclusion}
The problem of finding non-trivial bispectral situations remains a challenging one. In the case of \cite{DG}, it was only because one could painfully and explicitly solve the corresponding ad-conditions, namely 
$\operatorname {ad}_L^5(\Theta)=0$, that the first example was found. This opened the way to unexpected 
connections with the KdV equation and all the tools that this implies such as the Darboux
transformations. Nothing in the original ad-conditions made one suspect that such treasure was lying around.

One can only hope that being able to solve some of the new ad-conditions  of the kind considered here may be the key to unlocking new hidden gems, such an approach was used in \cite{GI03}.

The two extensions of the Reach ad-conditions given in this paper, to
Dunkl operators and to matrix-valued polynomials, have been treated
separately, and it is natural to ask what happens when they are combined.
The construction of \cite{B26} used in Examples~\ref{ex:laguerre-5}
and~\ref{ex:laguerre-7} starts from a pair of scalar weights $w_{1},w_{2}$
whose orthogonal polynomials are eigenfunctions of two second-order
operators with suitably matched eigenvalues. For the generalized Hermite
weights $w_{1}=|x|^{2\mu}e^{-x^{2}}$, $w_{2}=|x|^{-2\mu}e^{-x^{2}}$ of
Example~\ref{ex:gen-hermite}, with $|\mu|<\tfrac12$, this matching can be
arranged with the Dunkl operators $D_{\mu}$ and $D_{-\mu}$ in place of the
differential operators, and one obtains a candidate family of matrix-valued
generalized Hermite polynomials whose ad-conditions should display, at the
same time, the reflection $R$ of Section~\ref{sec:dunkl} and the matrix
structure of Section~\ref{sec:matrix}. A second natural direction is the
exceptional Dunkl case: a sequence $\{H^{\mu}_{n}\}_{n\in\mathbb{N}_{0}\setminus X}$
with a finite set $X$ of exceptional degrees should fit
Theorem~\ref{thm:BGZ-general}, with the coefficients $V_{j}(n)$ spectrally
realized by elements of $\mathcal{A}_{\mu}$ built from $D_{\mu}$ and $R$; since $R$
only distinguishes the two parities, the admissible sets $X$ are expected
to be more rigid than in the classical case. We plan to return to both
questions in a future publication.

\section*{Acknowledgments}

The research of I.~Bono Parisi was partially supported by SeCyT-UNC,
PIP 33620230100819CB, CONICET, PIP 1220150100356, and AUIP PMA2-2023-056-14.
The research of I.~Zurri\'an was partially supported by PID2021-124332NB-C21
(Ministerio de Ciencia e Innovaci\'on and FEDER Funds, European Union),
FQM-262 (Junta de Andaluc\'ia), and AUIP PMA2-2024-120-14.

\section*{Declaration on the use of AI tools}

In the preparation of this manuscript the authors used the AI assistant
Claude (Anthropic) and ChatGPT (OpenAI) to assist with the editing of parts of the text.
All mathematical results, proofs and conclusions are the authors' own, and
the authors take full responsibility for the content of the paper.

\end{document}